\documentclass{article}

\usepackage{lipsum}
\usepackage{subfigure}
\usepackage{amsfonts}
\usepackage{graphicx}
\usepackage{epstopdf}
\usepackage{amsfonts}
\usepackage{graphicx}
\usepackage{amsopn}
\usepackage{listings}
\usepackage{xcolor}
\usepackage{float}
\usepackage{marvosym}
\usepackage{booktabs}
\usepackage{multirow}

\usepackage{algorithm}
\usepackage{algorithmic}

\usepackage{graphicx}

\usepackage{lineno}
\modulolinenumbers[1]   

\usepackage{amsmath}
\usepackage{amsfonts}
\usepackage{amsbsy}
\usepackage{amscd}
\usepackage{amssymb}
\usepackage{latexsym}
\usepackage[mathscr]{eucal}
\usepackage{bm}

\usepackage{graphicx}
\usepackage{subfigure}
\usepackage{amsmath}	
\usepackage{amssymb}

\usepackage{listings}
\usepackage{xcolor}
\usepackage{float}
\usepackage{booktabs}
\usepackage{multirow}
\usepackage[noblocks]{authblk} %

\usepackage{algorithm}
\usepackage{algorithmic}
\floatname{algorithm}{Algorithm}

\newtheorem{proposition}{Proposition}[section]

\def\vecx{\mathbf{x}}
\def\vecy{\mathbf{y}}
\def\vecz{\mathbf{z}}
\def\vecr{\mathbf{r}}
\def\vect{\mathbf{t}}
\def\vecu{\mathbf{u}}
\def\vecv{\mathbf{v}}
\def\vecw{\mathbf{w}}

\def\matA{\mathbf{A}}
\def\matB{\mathbf{B}}
\def\matC{\mathbf{C}}
\def\matU{\mathbf{U}}
\def\matV{\mathbf{V}}
\def\matW{\mathbf{W}}
\def\matI{\mathbf{I}}

\def\df{\, \mathrm{d}}

\newcommand{\revise}[1]{{\color{black}#1}}
\newcommand{\rrevise}[1]{{\color{black}#1}}

\begin{document}




\title{A New Subspace Iteration Algorithm for Solving Generalized Eigenvalue Problems}


\author[a,b]{Biyi Wang}
\author[a,c]{Hengbin An \textsuperscript{\Letter} 
}
\author[d,e]{Hehu Xie 
}
\author[a,c]{Zeyao Mo 
}

\affil[a]{Institute of Applied Physics and Computational Mathematics, Beijing 100094, China}
\affil[b]{Graduate School of China Academy of Engineering Physics, Beijing 100088, China}
\affil[c]{CAEP Software Center for High Performance Numerical Simulation, Beijing 100088, China}
\affil[d]{ICMSEC, LSEC, NCMIS, Academy of Mathematics and Systems Science, Chinese Academy of Sciences, Beijing
100190, China}
\affil[e]{School of Mathematical Sciences, University of Chinese Academy of Sciences, Beijing 100049, China}

\renewcommand*{\Affilfont}{\small\it} 
\renewcommand\Authands{ and }

\date{}
\maketitle

\begin{abstract}
  It is needed to solve generalized eigenvalue problems (GEP) in many applications,
  such as the numerical simulation of vibration analysis, \revise{quantum mechanics, electronic structure},  etc.
  The subspace iteration is a kind of widely used algorithm to solve eigenvalue problems.
  To solve the generalized eigenvalue problem, one kind of subspace iteration method, Chebyshev-Davidson algorithm, is proposed recently.
  In Chebyshev-Davidson algorithm, the Chebyshev polynomial filter technique
  is incorporated in the subspace iteration~\cite{miao2020chebyshev}.
  In this paper, based on Chebyshev-Davidson algorithm,
  a new subspace iteration algorithm is constructed.
  In the new algorithm, the Chebyshev filter and
  inexact Rayleigh quotient iteration techniques are combined together
  to enlarge the subspace in the iteration.
  Numerical results of a vibration analysis problem show that
  the number of iteration and computing time
  of the proposed algorithm is much less than
  that of the Chebyshev-Davidson algorithm and some typical GEP solution algorithms.
  Furthermore, the new algorithm is more stable and reliable than
  the Chebyshev-Davidson algorithm in the numerical results.

\end{abstract}



\let\thefootnote\relax\footnotetext{
                This work was funded by National Natural Science Foundation of China (No. 12171045),
                Beijing Natural Science Foundation (No. Z200003), 
                and Science Challenge Project (No. TZ2016002).}




\section{Introduction}
\label{sect:introd}

Consider the following symmetric generalized eigenvalue problem (GEP)
\begin{eqnarray}
\label{eq:GEP}
\matA \vecx = \lambda \matB \vecx, \quad \text { with } \quad \| \vecx \|_{\matB} = 1,
\end{eqnarray}
where $\matA \in \mathbb{R}^{N \times N}$ and $\matB \in \mathbb{R}^{N \times N}$
are large, sparse and symmetric matrices with $\matB$ being positive definite,
$\| \vecx \|_{\matB}$ is the $\matB$-norm of the vector $\vecx$ which will be defined later.
In many applications, it is needed to solve GEP~\eqref{eq:GEP}, such as
the vibration analysis, quantum mechanics, electronic structure calculations, etc.
In~\eqref{eq:GEP}, the matrix $\matA$ is the stiffness matrix, and $\matB$ represents
the mass matrix. $\lambda$ is an eigenvalue of the matrix pencil $(\matA, \matB)$,
and the nonzero vector $\vecx$ is the corresponding eigenvector.
$(\lambda, \vecx)$ is called an eigenpair of GEP.
If $\matB$ is a diagonal matrix, then the generalized eigenvalue problem
reduces to the standard eigenvalue problem (SEP).
Usually, the smallest eigenvalues and the corresponding eigenvectors
play more important role in real applications than the largest ones.
We are interested in the computation for the smallest eigenvalues and
corresponding eigenvectors of a large sparse eigenvalue problem in this paper.

By now, many algorithms have been developed for solving eigenvalue problems~\cite{golub2000review,sorensen2002review,van2002review}.
These algorithms can be divided roughly into two classes:
the direct methods and iterative methods.
For large sparse eigenvalue problem, the iterative methods are preferable.
Three types of iterative methods are very popular.

The first one is based on the idea of the steepest descent method~\cite{steepestdescentmethod},
and the specific algorithms include
the conjugate gradient method~\cite{conjugategradient,projectedconjugategradient},
locally optimal block preconditioned conjugate gradient (LOBPCG) method~\cite{knyazev2001toward},
and generalized conjugate gradient eigensolver (GCGE)~\cite{li2021gcge}.

The second one is the Krylov subspace iteration methods,
including the Arnoldi method~\cite{arnoldi1951principle},
the Lanczos method~\cite{lanczos1950iteration},
and their important variants: the implicit restart Arnoldi method~\cite{sorensen1992implicit},
the rational Krylov method~\cite{rationalKrylov,ruhe1984rational} and
the Krylov-Schur method~\cite{stewart2002krylov}.
Also, some block versions Krylov subspace methods have been proposed.

The third one is Davidson-type subspace method which is more flexible
since there is no need to maintain the Krylov subspace structure in this type of methods.
In each iteration of Davidson-type methods, the subspace
is enlarged by adding an augmentation vector,
which is typically obtained by solving a correction equation. The performance of Davidson-type methods is strongly dependent on
the construction of the correction equation.
The original Davidson method~\cite{davidsorq1975theiterative} works well
only when the matrix is diagonally dominant.
In order to accelerate the convergence of the iteration, generalized Davidson (GD) method ~\cite{morgan1986generalizations}, Jacobi-Davidson (JD) method ~\cite{sleijpen2000jacobi} and
its variants~\cite{sleijpen1998jacobi,hochstenbach2003two,sleijpen1996jacobi}
are proposed with better correction equations.
Some efficient (preconditioned) linear equation solver
is needed to solve a correction equation.

In 2007, Zhou and Saad proposed a special kind of subspace iteration method
for solving standard eigenvalue problem by employing the polynomial filter technique
in the iteration. Specifically, they used the Chebyshev filter
in the subspace iteration and the obtained algorithm is called Chebyshev-Davidson method~\cite{zhou2007chebyshev}.
The performance of Chebyshev-Davidson method is excellent when
it is used to solve large sparse SEP.
Two years ago, Miao made some generalization of the Chebyshev-Davidson algorithm
so that it can be used to solve the generalized eigenvalue problem~\cite{miao2020chebyshev}.
The Chebyshev-Davidson method shows good effectiveness and robustness.
More importantly, in Chebyshev-Davidson method,
it is only needed to execute matrix-vector product,
and there is no need to construct or solve any correction equations.

On the other hand,
Zhou~\cite{zhou2006studies} revealed that the success of Davidson type methods
comes from the approximate Rayleigh quotient iteration (RQI) direction
implied in the solution of the correction equation.
The RQI is globally convergent for Hermitian matrices,
and it yields cubic convergence rate in some cases.
However, the RQI is not often used in practice because of
the high cost of the frequent factorizations.
To improve the efficiency of the RQI,
many researchers proposed to solve the equations
inexactly by some Krylov iterative methods in RQI,
such as FGMRES~\cite{szyld2011efficient},
MINRES~\cite{jia2012convergence}, and CRM~\cite{simoncini2002inexact}.
At the same time, some efficient preconditioners were also constructed.
However, to make the algorithms work efficiently,
either a good approximate eigenpair should be provided,
or the shifted matrix should be decomposed to construct a preconditioner.
In real applications, it is usually difficult to give a good
approximate eigenpair, and for solving large scale problem,
it is too expensive to decompose the shifted matrix.
Furthermore, they only considered to compute one
eigenvalue and the corresponding eigenvector.

In this paper, by combining the \revise{Chebyshev filtering  technique and inexact Rayleigh quotient iteration}, 
a new subspace iteration method is proposed to solve GEP,
and the obtained algorithms is called the Chebyshev-RQI subspace
(CRS) iteration method.
Numerical results show that the number of iteration of Chebyshev-RQI
subspace method is much less than the Chebyshev-Davidson and some
typical GEP solution algorithms.
Particularly, CRS shows much better performance than CD
in terms of computing time in the case of parallel computing,
which is indispensable to large scale numerical simulations.

This paper is organized as follows.
In Sect.~\ref{sect:fundamentals}, the polynomial filtering technique
for solving symmetric generalized eigenvalue problems 
and Rayleigh quotient iteration are introduced.
The Chebyshev-Davidson method for solving generalized eigenvalue problem
is given in Sect.~\ref{sect:Cheby-davidson}.
The Chebyshev-RQI subspace method is proposed in Sect.~\ref{sect:CRS method}.
In Sect.~\ref{sect:numer-result}, some numerical results are given to
show the effectiveness of the proposed algorithm.
Finally, some conclusions and remarks are given in Sect.~\ref{sect:conclusion}.

\paragraph{Notations and basic assumptions}
Throughout this paper,
$\matA \in \mathbb{R}^{N \times N}$ and $\matB \in \mathbb{R}^{N \times N}$
represent large, sparse and symmetric matrices, and the matrix $\matB$ is positive definite.
For vectors $\vecu$ and $\vecv$ in the real vector space $\mathbb{R}^{N}$,
$\langle \vecu, \vecv \rangle$ represents the inner product,
and the $\matB$-inner product with respect to the symmetric and
positive definite matrix $\matB$ is defined by
$
\langle \vecu, \vecv \rangle_{\matB} =
\langle \vecu, \matB \vecv \rangle =
\vecv^{\top} \matB \vecu.
$
If $\langle \vecu, \vecv \rangle_{\matB} = 0$, then we say that vector $\vecu$ is $\matB$-orthogonal to vector $\vecv$,
denoted by $\vecu \perp_{\matB} \vecv$. Here and in the subsequent,
the superscript $(\cdot)^{\top}$ denotes the transpose of
either a matrix or a vector. Thus, the corresponding $\matB$-norm of
a vector $\vecv$ is defined by $\| \vecv \|_{\matB} = \sqrt{\langle \vecv, \vecv \rangle_{\matB}}$.
In particular, when $\matB$ is the identity matrix,
the $\matB$-norm of a vector $\vecv$ reduces to the standard Euclidean norm $\| \vecv \|$.
Denote by $\mathcal{P}_{m}$ the set of all polynomials with degree $m$.
The spectrum of the matrix pencil $(\matA, \matB)$ is denoted by $\left\{ \lambda_{i} \right\}_{i=1}^{N}$
with ascending order, i.e., $\lambda_{1} < \lambda_{2} \leq \cdots \leq \lambda_{N}$,
and denote by $\left\{ \vecx_{i} \right\}_{i=1}^{N}$ the corresponding eigenvectors with
$\| \vecx_i \|_{\matB} = 1$, $i = 1, 2, \ldots, N$.
For a nonzero vector $\vecz$, the value $\theta_{\vecz} = \frac{\vecz^{\top} \matA \vecz}{\vecz^{\top} \matB \vecz}$ is called the
generalized Rayleigh quotient of $\vecz$ associated with matrix pencil $(\matA, \matB)$.

\section{Chebyshev Filtering and Rayleigh Quotient Iteration}
\label{sect:fundamentals}

The Chebyshev polynomial is a kind of important polynomial
that \revise{is widely used in solving} large scale linear equations and
eigenvalue problems. The Rayleigh quotient iteration is a kind of elementary method
for solving eigenvalue problems. The Chebyshev polynomial filtering technique
and Rayleigh quotient iteration are the basis \revise{for designing the new} algorithm in this paper.

\subsection{Chebyshev Filtering}
\label{subsect:chebyshev-filter}

The Chebyshev filtering technique for solving eigenvalue problem is based on
the Chebyshev polynomials.
The real Chebyshev polynomials of the first kind are defined as follows~\cite{saad2011numerical,zhou2007chebyshev}
\begin{eqnarray*}
\label{eqn:Chebyshev-polynomial}
C_{m}(t)=
\begin{cases}
\cos \left(m \cos ^{-1}(t)\right),   & -1 \leq t \leq 1, \\
\cosh \left(m \cosh ^{-1}(t)\right), & |t|>1 ,
\end{cases}
\end{eqnarray*}
where $m=0, 1, \ldots$.
Note that $C_{0}(t)=1$, $C_{1}(t)=t$, and $C_{m}(t)$ can be computed straightforwardly
by the three-term recurrence \revise{formula}
\begin{eqnarray*}
C_{m+1}(t)= 2 t C_{m}(t) - C_{m-1}(t), \quad m = 1, 2, \ldots.
\end{eqnarray*}
Besides, Chebyshev polynomials have another appealing character
with the rapid growth outside the interval [-1,1].
For more details and \revise{properties about Chebyshev polynomials,
please refer to}~\cite{saad2011numerical,zhou2007chebyshev}.

The Chebyshev polynomial filtering technique has been explored thoroughly
for standard eigenvalue problem~\cite{zhou2007chebyshev}.
Assume that $(\mu, \vecw)$ is an eigenpair of the matrix $\matA$,
and $p(t) \in \mathcal{P}_m$ is a polynomial,
where $m$ is an prescribed integer. \revise{Then} it is easy to know that
\begin{eqnarray}
\label{eqn:eigenvalue-polynomial}
p(\matA) \vecw = p(\mu) \vecw.
\end{eqnarray}
This is the basis that the Chebyshev polynomial can be used to
filter \revise{out the components in unwanted eigenspaces for the} standard eigenvalue problem
$\matA \vecx = \lambda \vecx$~\cite{saad2011numerical}.
For generalized eigenvalue problems~\eqref{eq:GEP},
however, the Chebyshev polynomial filtering technique can not be used directly
because it is not easy to obtain similar relationship as~\eqref{eqn:eigenvalue-polynomial}
for generalized eigenvalue problem.
\revise{Recently}, Miao generalized the Chebyshev filtered iteration
to the case of the matrix pencil $(\matA, \matB)$ for generalized eigenvalue problem~\cite{miao2020chebyshev}.
The generalization is based on the following proposition.

\begin{proposition}\label{Proposition_1}
Assume that $\left(\lambda_{1}, \vecx_{1}\right)$ with
$\left\| \vecx_{1} \right\|_{\matB} = 1$ is an eigenpair
of the matrix pencil $(\matA, \matB)$. Then $\vecx_{1}$ is also an eigenvector of the symmetric matrix $\matA - \lambda_{1} \matB$
associated with the zero eigenvalue.
\end{proposition}

By Proposition~\ref{Proposition_1}, if $\lambda_1$ is available,
then we can use the approximate eigenvector associated
with the zero eigenvalue of the symmetric matrix $\matA - \lambda_{1} \matB$
to approximate $\vecx_{1}$.
Thus we can convert a generalized eigenvalue problem to a standard eigenvalue problem
and then we could consider the application of Chebyshev filter on the symmetric matrix
$\matA - \lambda_{1} \matB$.

Since $\matA - \lambda_{1} \matB$ is a real symmetric semi-positive matrix,
and $\vecx_1$ is the first eigenvector of this matrix,
it is easy to know that there exists an orthogonal matrix $\matU$ with
\begin{eqnarray}
\label{eqn:matrix-U-form}
\matU = \left[ \vecu_{1}, \vecu_{2}, \ldots, \vecu_{N} \right],
\quad
\text{where}
\quad
\vecu_{1} = \frac{\vecx_{1}}{\left\| \vecx_{1} \right\|},
\end{eqnarray}
such that
\begin{eqnarray}
\label{eqn:decompos-A-theta-B}
\matA - \lambda_{1} \matB = \matU \mathbf{\Sigma} \matU^{\top},
\end{eqnarray}
where $\mathbf{\Sigma}$ is a diagonal matrix with following form
\begin{eqnarray*}
\mathbf{\Sigma} = \operatorname{diag} \left( \sigma_{1}, \sigma_{2}, \ldots, \sigma_{N} \right),
\quad
0 = \sigma_{1} < \sigma_{2} \leq \cdots \leq \sigma_{N}.
\end{eqnarray*}
Since $\{ \vecx_j \}_{j=1}^N$ is an orthogonal basis for $\mathbb{R}^N$,
therefore, $\vecu_i$ can be expressed by
\begin{eqnarray}
\label{eqn:u-express-by-x}
\vecu_i = \sum_{j=1}^N \kappa_{i,j} \vecx_j,
\quad
\kappa_{i,j} \in \mathbb{R},
\quad
i = 1, 2, \ldots, N.
\end{eqnarray}
In particular, by~\eqref{eqn:matrix-U-form},
\begin{eqnarray}
\label{eqn:kappa-1j}
\kappa_{1,1} = \frac{1}{\| \vecx_1 \|},
\quad
\kappa_{1,j} = 0,
\quad
j = 2, \ldots, N.
\end{eqnarray}

Next assume that there exists an approximate eigenvector $\tilde{\vecx}$ for $\vecx_1$
(for example, an approximation can be obtained in the iteration of some subspace method).
Since $\tilde{\vecx} \approx \vecx_1$, therefore, by~\eqref{eqn:matrix-U-form}, we have
\begin{eqnarray}
\label{eqn:u1-approx-x1}
\left\langle \vecu_1, \tilde{\vecx} \right\rangle
=
\left\langle \frac{\vecx_1}{\left\| \vecx_1 \right\|}, \tilde{\vecx} \right\rangle
\approx \| \vecx_1 \|.
\end{eqnarray}

Let $p(t)$ be a polynomial.
By this polynomial and $\tilde{\vecx}$, we hope to construct a better approximate vector
$\tilde{\vecz}$ for $\vecx_1$ ($\tilde{\vecz}$ is called a polynomial filtered vector).
Let $\tilde{\vecz}$ defined by
\begin{eqnarray*}
\label{eqn:filter-vec-z}
\tilde{\vecz} = p \left( \matA - \lambda_1 \matB \right) \tilde{\vecx}.
\end{eqnarray*}
Then by~\eqref{eqn:decompos-A-theta-B},
\eqref{eqn:u-express-by-x},
\eqref{eqn:kappa-1j},
\eqref{eqn:u1-approx-x1}
and the orthogonality of the matrix $\matU$, we have
\begin{eqnarray*}
\tilde{\vecz}
&=& p \left( \matA - \lambda_1 \matB \right) \tilde{\vecx} 
= p \left( \matU \mathbf{\Sigma} \matU^{\top} \right) \tilde{\vecx} 
= \matU p \left( \mathbf{\Sigma} \right) \matU^{\top} \tilde{\vecx} 
= \sum_{i=1}^N p \left( \sigma_{i} \right)
    \left\langle \vecu_i, \tilde{\vecx} \right\rangle \vecu_i \\
&=& \sum_{i=1}^N p \left( \sigma_{i} \right)
    \left\langle \vecu_i, \tilde{\vecx} \right\rangle
    \sum_{j=1}^N \kappa_{i,j} \vecx_j 
= \sum_{j=1}^N
    \left(
    \sum_{i=1}^N p \left( \sigma_{i} \right) \left\langle \vecu_i, \tilde{\vecx} \right\rangle \kappa_{i,j}
    \right) \vecx_j \\
&=& \left(
    \sum_{i=1}^N p \left( \sigma_{i} \right) \left\langle \vecu_i, \tilde{\vecx} \right\rangle \kappa_{i,1}
    \right) \vecx_1 +
    \sum_{j=2}^N
    \left(
    \sum_{i=1}^N p \left( \sigma_{i} \right) \left\langle \vecu_i, \tilde{\vecx} \right\rangle \kappa_{i,j}
    \right) \vecx_j \\
&=&
    \left(
    p \left( \sigma_{1} \right) \left\langle \vecu_1, \tilde{\vecx} \right\rangle \kappa_{1,1} +
    \sum_{i=2}^N p \left( \sigma_{i} \right) \left\langle \vecu_i, \tilde{\vecx} \right\rangle \kappa_{i,1}
    \right) \vecx_1 +
    \sum_{j=2}^N
    \left(
    \sum_{i=2}^N p \left( \sigma_{i} \right) \left\langle \vecu_i, \tilde{\vecx} \right\rangle \kappa_{i,j}
    \right) \vecx_j \\
&\approx&
    \left(
    p \left( \sigma_{1} \right) \| \vecx_1 \| \frac{1}{\| \vecx_1 \|} +
    \sum_{i=2}^N p \left( \sigma_{i} \right) \left\langle \vecu_i, \tilde{\vecx} \right\rangle \kappa_{i,1}
    \right) \vecx_1 +
    \sum_{j=2}^N
    \left(
    \sum_{i=2}^N p \left( \sigma_{i} \right) \left\langle \vecu_i, \tilde{\vecx} \right\rangle \kappa_{i,j}
    \right) \vecx_j \\
&=&
    \left(
    p \left( \sigma_{1} \right) +
    \sum_{i=2}^N p \left( \sigma_{i} \right) \left\langle \vecu_i, \tilde{\vecx} \right\rangle \kappa_{i,1}
    \right) \vecx_1 +
    \sum_{j=2}^N
    \left(
    \sum_{i=2}^N p \left( \sigma_{i} \right) \left\langle \vecu_i, \tilde{\vecx} \right\rangle \kappa_{i,j}
    \right) \vecx_j.
\end{eqnarray*}
In order to make the \revise{component of   $\tilde{\vecz}$ in the direction $\vecx_1$ 
as relatively large as possible}, we may choose a polynomial $p(t)$ such that
\begin{eqnarray*}
\left| p \left(\sigma_{1} \right) \right| \gg \left| p \left(\sigma_{i} \right) \right|,
\quad
i = 2, 3, \ldots, N.
\end{eqnarray*}
Similar to the derivation for the standard symmetric eigenvalue problems~\cite{saad2011numerical},
we can also define a min-max problem:
find $p(t) \in \mathcal{P}_{m}$ such that $p (\sigma_1) = 1$ and
\begin{eqnarray}
\label{eqn:min-max-prob}
\max _{t \in[a, b]}|p(t)|
=
\min _{q \in \mathcal{P}_{m} \atop q (\sigma_1)=1}
\max _{t \in[a, b]}|q(t)|,
\end{eqnarray}
where $[a, b]$ is an interval containing the eigenvalues
$\left\{ \sigma_{i} \right\}_{i=2}^{N}$ while excluding $\sigma_{1}$.
The polynomial $p(t)$ satisfying~\eqref{eqn:min-max-prob} is
a filter with degree $m$,
which can be prescribed as
\begin{eqnarray}
\label{scaledCheby}
p(t) = \frac{C_m \left( 1 + 2 \frac{t-b}{b-a} \right)}
            {C_m \left( 1 + 2 \frac{\sigma_{1} - b}{b-a} \right)},
\end{eqnarray}
where $C_{m}(t)$ is the Chebyshev polynomial of the first kind with degree $m$.

The above discussion is based on the assumption that the eigenvalue $\lambda_1$ is known.
Although $\lambda_{1}$ is unknown in the actual calculation,
some approximate eigenvalue can be given in the iteration of a specific method.
For example, the Rayleigh quotient at $\tilde{\vecx}$ can be used to approximate  the eigenvalue $\lambda_1$
\begin{eqnarray*}
\tilde{\theta} = \frac{ \tilde{\vecx}^{\top} \matA \tilde{\vecx}}
                      { \tilde{\vecx}^{\top} \matB \tilde{\vecx}}.
\end{eqnarray*}

In the following, assume that
an approximate eigenpair $(\tilde{\theta}, \tilde{\vecx})$ is available.
Let the eigenvalues of the matrix $\matA - \tilde{\theta} \matB$ are
\begin{eqnarray*}
\tilde{\sigma}_1 < \tilde{\sigma}_2 \le \cdots \le \tilde{\sigma}_N.
\end{eqnarray*}

In actual computation, the polynomial filtered process
$\tilde{\vecz} = p(\matA - \tilde{\theta} \matB) \tilde{\vecx}$
with the scaled Chebyshev polynomial
can be implemented once $\tilde{\sigma}_1$ is available,
and an interval $[a, b]$ can be determined such that
$\left\{ \tilde{\sigma}_2, \ldots, \tilde{\sigma}_N \right\} \subset [a, b]$.
Specifically, if $\tilde{\sigma}_1$, $a$, and $b$ are determined,
then the polynomial~\eqref{scaledCheby} is determined by setting $\sigma_1$
as $\tilde{\sigma}_1$, and
the polynomial filtered process can be implemented economically
with the three-term recurrence relation that can be described algorithmically as follows~\cite{saad2011numerical,zhou2007chebyshev}.

\begin{algorithm}[!htbp]
	\caption{Chebyshev Filtered Iteration}
	\label{alg:Chebyshev Filtered}
	\begin{algorithmic}[1]
		\REQUIRE $\matC$, $\vecx$, $m$, $a$, $b$, $\tilde{\sigma}_1$. \\
		\ENSURE $\vecz_{m}$ \\
		\STATE Set $\vecz_0 = \vecx$.
		\STATE Compute
		$$
		\mu = \frac{b+a}{2}, \quad
		\nu = \frac{b-a}{2}, \quad
		\gamma_1 = \frac{\nu}{\tilde{\sigma}_1 - \mu}, \quad
		\text{and} \quad
		\vecz_1 = \frac{\gamma_1}{\nu}(\matC \vecx - \mu \vecx).
		$$
		\FOR{$i=1,2,\dots, m-1$}
		\STATE Set $\gamma_{i+1} = \frac{1}{2 / \gamma_{1} - \gamma_{i}}$;
		\STATE Compute
		$$
		\vecz_{i+1} = 2 \frac{\gamma_{i+1}}{\nu}
		\left(\matC \vecz_{i} - \mu \vecz_{i} \right)
		- \gamma_{i} \gamma_{i+1} \vecz_{i-1}.
		$$
		\ENDFOR
	\end{algorithmic}
\end{algorithm}

\subsection{Rayleigh Quotient Iteration}
\label{subsect:RQI}

The Rayleigh quotient iteration is \revise{an important accelerating technique}
for solving eigenvalue problem.
Assume that $\matC$ is a normal matrix,
the Rayleigh quotient iteration for standard eigenvalue problem $\matC \vecx = \lambda \vecx$
can be described by \revise{Algorithm \ref{alg:RQI}}.

\begin{algorithm}[!htbp]
	\caption{Rayleigh Quotient Iteration}
	\label{alg:RQI}
	\begin{algorithmic}[1]
		\REQUIRE Matrix $\matC$, \revise{initial vector} $\vecv_0$,
                 convergence tolerance $\varepsilon_{\text{\tiny{RQI}}}$. \\
		\ENSURE $\vecv$. \\
  		\STATE $\vecv_0 = \vecv_0 / \|\vecv_0\|$.
		\FOR{$i = 1,2,\dots$}
		\STATE $\tau_i = \left\langle \matC \vecv_{i-1}, \vecv_{i-1} \right\rangle$.
		\STATE $\hat{\vecv}_i = \left( \matC - \tau_i \matI \right)^{-1} \vecv_{i-1}$.
		\STATE $\vecv_i = \hat{\vecv}_i / \|\hat{\vecv}_i\|$.
        \IF{$\Vert \matC \vecv_i - \tau_i \vecv_i \Vert < \varepsilon_{\text{\tiny{RQI}}}$}
		\STATE $\vecv = \vecv_i$;
		\STATE break;
		\ENDIF
		\ENDFOR
	\end{algorithmic}
\end{algorithm}

The scalar $\tau_i$ on Line 3 in Algorithm~\ref{alg:RQI} is the Rayleigh quotient.
The initial vector $\vecv_0$ has important influence on the behaviour of the algorithm.
The iteration sequences produced by this algorithm usually converge to
the eigenvector which is close to the initial vector.
To compute $\hat{\vecv}_i$ on Line~4, it is needed to solve
a linear system $\left( \matC - \tau_i \matI \right) \hat{\vecv}_i = \vecv_{i-1}$ exactly, \revise{which is very expensive for large scale problems}.
When an iterative method, such as a Krylov subspace method,
is used to solve this linear system, the deduced  algorithm is called
inexact Rayleigh quotient iteration (IRQI) method, which is described by Algorithm~\ref{alg:IRQI}.
This algorithm is an inner-outer iteration process with the outer iteration be the Rayleigh quotient iteration and the inner iteration for solving the shifted linear system.

\begin{algorithm}[!htbp]
	\caption{Inexact Rayleigh Quotient Iteration}
	\label{alg:IRQI}
	\begin{algorithmic}[1]
		\REQUIRE Matrix $\matC$, \revise{initial vector}  $\vecv_0$,
                 maximal iterate number $\rm It_{\text{max-\tiny{RQI}}}$,
                 convergence tolerance $\varepsilon_{\text{\tiny{RQI}}}$,
                 and maximal iteration number $\rm It_{\text{max-linear}}$ for solving linear equation.
		\ENSURE $\vecv$. \\
  		\STATE $\vecv_0 = \vecv_0 / \| \vecv_0 \|$.		
		\FOR{$i=1,\dots, It_{\text{max-\tiny{RQI}}}$}
		\STATE $\tau_i = \left\langle \matC \vecv_{i-1}, \vecv_{i-1} \right\rangle$.
		\STATE Solve $\left( \matC - \tau_i \matI \right) \hat{\vecv}_i = \vecv_{i-1}$
               by some iterative method with maximal $\rm It_{\text{max-linear}}$ \revise{iteration steps}.
		\STATE $\vecv_i = \hat{\vecv}_i / \|\hat{\vecv}_i\|$.
        \IF{$\Vert \matC \vecv_i - \tau_i \vecv_i \Vert < \varepsilon_{\text{\tiny{RQI}}}$}
		\STATE $\vecv = \vecv_i$;
		\STATE break;
		\ENDIF
		\ENDFOR
	\end{algorithmic}
\end{algorithm}

\section{Chebyshev Davidson Method}
\label{sect:Cheby-davidson}

In this section, we introduce briefly the Chebyshev-Davidson method
for GEP.
The Chebyshev-Davidson method is a kind of subspace iteration algorithm,
which extracts approximate eigenpairs by the Rayleigh-Ritz method~\cite{parlett1998symmetric}
on a sequence of gradually enlarging subspaces.
The most important component affecting the performance of the algorithm
is the construction of a proper polynomial filter in each iteration.
Specifically, the choice of a proper filter interval and a polynomial order.
Miao~\cite{miao2020chebyshev} gives \revise{a simple and  practical advice} in his paper without estimating the upper bound
of the eigenvalues.

The Chebyshev-Davidson method~\cite{miao2020chebyshev} is described by Algorithm~\ref{alg:Chebyshev-Davidson}.
This algorithm compute altogether $NEV$ eigenpairs.
To compute each eigenpair, a subspace $\mathbf{\mathcal{V}}$ is constructed by an iteration process
(with the iteration index $k$).
In the Chebyshev filtered process of this algorithm,
the smallest eigenvalue $\tilde{\sigma}_{1}$ of the matrix $\matA - \theta^{(k)} \matB$
and the interval $[a, b]$ containing $\left\{ \tilde{\sigma}_{i} \right\}_{i=2}^{N}$
should be prescribed in advance,
where $\theta^{(k)}$ is the approximate eigenvalue of the matrix pair $(\matA, \matB)$ in the $k$-th iteration,
and $\left\{ \tilde{\sigma}_{i} \right\}_{i=1}^{N}$ are the eigenvalues of the matrix $\matA - \theta^{(k)} \matB$.
For this purpose, assume that after the $k$-th iteration, a subspace $\mathbf{\mathcal{V}}_{k+1}$ is obtained,
and let $\matV_{k+1} \in \mathbb{R}^{N \times {(k+1)}}$ be a basis of this subspace.
Then, we may
let $\tilde{\sigma}_{1}, a$ and $b$ be the smallest,
the second smallest and the largest eigenvalues of the projected matrix
\begin{eqnarray*}
\matV_{k+1}^{\top} \matC \matV_{k+1}
: =
\matV_{k+1}^{\top} \left( \matA - \theta^{(k)} \matB \right) \matV_{k+1}
=
\matV_{k+1}^{\top} \matA \matV_{k+1} - \theta^{(k)} \matV_{k+1}^{\top}
\matB \matV_{k+1}.
\end{eqnarray*}

\begin{algorithm}[!htbp]
	\caption{Chebyshev-Davidson Method for GEP}
	\label{alg:Chebyshev-Davidson}
	\begin{algorithmic}[1]
		\REQUIRE $\matA$, $\matB$, $NEV$, $m$, \rrevise{ $\rm dim_{\max}$}, $\rm It_{\max}$, $\varepsilon_{\text{\tiny{CD}}}$, $\vecv_0$.\\ %
		\ENSURE $\mathbf{\Lambda}$, $\matW$, $\rm It_{\text{total}}$.\\
		\STATE $\vecx = \vecv_0, \rm It_{\text{total}} = 0, \matW=[\ ], \mathbf{\Lambda} = [\ ]$.
		\FOR{$n = 1$ to \textit{NEV}}
		\STATE Let $k = 1$.
		\STATE Let $a = 0$, $b = 0$, and $\tilde{\sigma}_1 = 0$.
		\STATE Let $\matV_k = \{ \vecx \}$, $\theta^{(k)}_n = \frac{\vecx^{\top} \matA \vecx}{\vecx^{\top} \matB \vecx}$,
                     \rrevise{$\vecr = \frac{(\matA - \theta^{(k)}_n \matB) \vecx}{\theta^{(k)}_n \Vert \vecx \Vert}$.}
		\FOR{$k = 1$ to $\rm It_{\text{max}}$}
		\IF{$ \revise{ \Vert \vecr \Vert} < \varepsilon_{\text{\tiny{CD}}}$}
		\STATE $\rm It_{\text{total}} = k + \rm It_{\text{total}}$,
               $\matW = [\matW, \vecx]$, $\mathbf{\Lambda} = [\mathbf{\Lambda}, \theta^{(k)}_n]$.
		\STATE break.
		\ENDIF
		\IF{$\text{\rm dim}(\matV_k) \geq {\rm dim}_{\max}$}
		\STATE $\texttt{RESTART}: \matV_k = [\vecx], \theta^{(k)}_n = \frac{\vecx^{\top} \matA \vecx}{\vecx^{\top} \matB \vecx}$.
		\ENDIF
		\STATE Update $\matC = \matA - \theta^{(k)}_n \matB$.
		\IF{$k = 1$}
		   \STATE $\vecz = \matC \vecx$.
		\ELSE
		   \STATE $\vecz = {\texttt{ChebyshevFilter}}(\matC, \vecx, m, a, b, \tilde{\sigma}_1)$.
		\ENDIF
		\STATE Orthonormalize $\vecz$ against $\matW$ and $\matV_k$, respectively,
		       and update $\matV_{k+1} = [\matV_k, \vecz]$ and
               $\matV_{k+1}^{\top} \matC \matV_{k+1} = \matV_{k+1}^{\top} \matA \matV_{k+1} - \theta^{(k)}_n \matV_{k+1}^{\top} \matB \matV_{k+1}$.
		\STATE Compute the eigenvalues of $\matV_{k+1}^{\top} \matC \matV_{k+1}$,
               and determine the filter parameters
		       $a$, $b$, and $\tilde{\sigma}_1$ by
               the second smallest, the largest, and the smallest eigenvalues.
		\STATE Compute the eigenpairs $(\mu_i, \vecy _{i})$ of the matrix pair
               $(\matV_{k+1}^{\top} \matA \matV_{k+1}, \matV_{k+1}^{\top} \matB \matV_{k+1})$
		       in ascending sort order, $i=1, 2, \dots, k+1$.
		\STATE Update $\vecx = \matV_{k+1} {\vecy_1}$ and $\theta^{(k+1)}_n = \mu_1$.
        \STATE \revise{Update $\vecr = \frac{(\matA - \theta^{(k+1)}_n \matB) \vecx}{\theta^{(k+1)}_n \Vert \vecx \Vert}$.}
		\ENDFOR
		\STATE $\matV^{n} = \matV_{k}$, $\vecx = \matV^{n} {\vecy_2}$.
		\ENDFOR
	\end{algorithmic}
\end{algorithm}

In Algorithm~\ref{alg:Chebyshev-Davidson},
$NEV$ denotes the number of eigenpairs to be computed.
All these $NEV$ eigenpairs are computed one by one,
and the smallest eigenvalue is firstly obtained.
The maximal dimension of the subspace $\mathbf{\mathcal{V}}$ is limited by the parameter $dim_{\max}$.
The convergence tolerance for the subspace iteration is $\varepsilon_{\text{\tiny{CD}}}$,
and the maximal number of iteration is limited by the parameter $\rm It_{\max}$.
$\vecv_0$ is the initial vector of the iteration,
and usually set as a random unit vector.
$m$ is the degree of Chebyshev polynomial.
The converged eigenvalues are saved in $\mathbf{\Lambda}$,
and the converged eigenvectors are saved in $\matW$.
The total number of iteration is denoted by $It_{\text{total}}$.
On Line~18, {\texttt{ChebyshevFilter}} \revise{is defined by} Algorithm~\ref{alg:Chebyshev Filtered}.

It should be pointed out that Line~22 implies
\begin{eqnarray*}
\matV_{k+1}^{\top} \matA \matV_{k+1} \vecy_i =
\mu_i \matV_{k+1}^{\top} \matB \matV_{k+1} \vecy_i,
\end{eqnarray*}
and
\begin{eqnarray*}
\mu_1 \leq \mu_2 \leq \ldots \leq \mu_{k+1}.
\end{eqnarray*}
On Line~26, $\matV^{n}$ is the basis of the subspace for computing the the $n$-th eigenpair,
and $\vecy_2$ is the eigenvector corresponding to
the second smallest eigenvalue of the matrix pair
$(({\matV^{n}})^{\top} \matA \matV^n, ({\matV^{n}})^{\top} \matB \matV^n)$.
To compute the $(n+1)$-th eigenpair, it prefers to
use $\matV^{n} \vecy_2$ to start the iteration process
because it is a good approximation of $\vecx_{n+1}$
after finishing computing the $n$-th eigenpair.

\section{Chebyshev-RQI Subspace Method}
\label{sect:CRS method}

One of the advantages of the Chebyshev-Davidson method
is that the implementation of the algorithm
only \revise{includes} the matrix-vector products.
At the same time, it is known that the
RQI method has the advantage of cubic convergence rate in some cases.
Based on this observation, we propose a new subspace method
by combining the Chebyshev filtering and RQI.
Specifically, we use Chebyshev filtering and RQI at the same time
to produce the augmentation vectors in the subspace iteration,
and the obtained algorithm is called Chebyshev-RQI subspace (CRS) iteration method.
Compared with Chebyshev-Davidson method, a new augmentation vector $\vect^{(k)}$
is introduced after adding the filtered vector $\vecz^{(k)}$ to the subspace
in each iteration of the Chebyshev-RQI subspace method.

The subspace of CRS algorithm is defined as
\begin{eqnarray*}
&&\mathcal{V}_{1} := \left\{\vecx^{(1)}\right\}, \\
&&\mathcal{V}_{k+1} := \text{span} \left\{ \mathcal{V}_k,
        p_k \left( \matA - \theta^{(k)} \matB \right) \vecx^{(k)}, \vect^{(k)} \right\}, \quad k = 1, 2, 3, \ldots,
\end{eqnarray*}
where $\theta^{(j)}$ is a prescribed proper shift, and
$p_j(t)$ is the prescribed polynomial filter, $j=$ $1,2, \ldots, k$.
$\vecx^{(1)}$ is the initial guess of the eigenvector,
and the sequences $\{ \vecx^{(j)} \}$ will approximate the desired eigenvector.
The vector $\vect^{(k)}$ is used as an augmentation vector for constructing the subspace $\mathcal{V}_{k+1}$.

Next we discuss the construction of the augmentation vector $\vect^{(k)}$ at $k$-th iteration.
For simplicity, we will only discuss the case of calculating the first smallest eigenpair $(\lambda_1, \vecx_1)$.
And for other eigenpairs, the discussion would be similar.

Recalling that
\begin{eqnarray*}
    \matA \vecx_1  = \lambda_1 \matB \vecx_1 \quad \Longleftrightarrow
    \quad (\matA - \lambda_1 \matB) \vecx_1 = 0 \cdot \vecx_1,
\end{eqnarray*}
\revise{which} means computing the desired eigenvector $\vecx_{1}$ of
the matrix pencil $(\matA, \matB)$ is equivalent to
compute the eigenvector associated with the zero eigenvalue
of the symmetric matrix $\matA - \lambda_{1} \matB$.
Assume that $\theta^{(k)}\approx \lambda_1$ and $\vecx^{(k)} \approx \vecx_1$
at $k$-th iteration,
so the matrix $\matC^{(k)} = \matA - \theta^{(k)} \matB$ has an eigenvalue
$\xi^{(k)} \approx 0$ and the corresponding eigenvector $\vect^{(k)}$.

When a good approximate vector is available for the desired eigenpair,
the RQI algorithm could achieve fast convergence.
And fortunately several Chebyshev filter \revise{steps}
could provide good one close to the desired eigenvector. 
\revise{So} we try to use RQI to compute $(\xi^{(k)}, \vect^{(k)})$.
Let the iteration sequences be $\left \{ \xi^{(k)}_i , \vect^{(k)}_i \right \}, i = 0, 1, 2, \ldots$,
then by referring to Algorithm~\ref{alg:RQI},
we can get the following iterative formulation
\begin{eqnarray*}
\left\{
\begin{array}{rcl}
\vect^{(k)}_{i+1} &=& \left[ (\matA - \theta^{(k)} \matB) - \tau_i \matI  \right] ^{-1} \vect^{(k)}_i, \quad i = 0, 1, 2, \ldots, \vspace{3pt} \\
\tau_i            &=& \left\langle  (\matA - \theta^{(k)} \matB) \vect^{(k)}_i, \vect^{(k)}_i \right\rangle, \quad i = 0, 1, 2, \ldots, \vspace{3pt} \\
\vect^{(k)}_0     &=& \vecx^{(k)}.
\end{array}
\right.
\end{eqnarray*}
Here we consider the single-step RQI, then we have
\begin{eqnarray*}
\label{single-step RQI}
\vect^{(k)} = \left[ \left( \matA - \theta^{(k)} \matB \right) - \tau_0 \matI  \right] ^{-1} \vecx^{(k)}.
\end{eqnarray*}
Besides, we know that $\tau_0 \approx \xi^{(k)} \approx 0$,
thus we consider acquiring the augmentation vector $\vect^{(k)}$ by
\begin{eqnarray}
\label{single-step IRQI}
	\vect^{(k)} = \left( \matA - \theta^{(k)} \matB \right)^{-1} \vecx^{(k)}.
\end{eqnarray}
For large scale problems, it will be very expensive to solve~\eqref{single-step IRQI} exactly by using direct method. Therefore, some iterative methods, particularly, the Krylov methods
with or without preconditioner,
can be used to solve this equation inexactly by setting a fixed number of inner iteration
or a low tolerance. This is called the inexact Rayleigh quotient iteration (IRQI) corresponding to Algorithm~\ref{alg:IRQI} with only one-step iteration.

It should be pointed out that
Zhou~\cite{zhou2006studies} had revealed that the success of
Davidson-type methods comes from the approximate Rayleigh quotient iteration direction
implied in the solution of the Davidson correction equation.
In other words, the direction of the approximate Rayleigh quotient is
the essential ingredient of the Davidson-type methods.
\rrevise{In Tables~\ref{tab:scale1}--\ref{tab:scale3}, we can see the the performance of CRS method will be significantly improved 
by adding the augmented vector $\vect^{(k)}$ than that of CD.}

The \revise{detailed} description of the Chebyshev-RQI subspace method is given by Algorithm~\ref{alg:Chebyshev-RQI Subspace}.

\begin{algorithm}[htbp]
	\caption{Chebyshev-RQI Subspace(CRS) Method for GEP}
	\label{alg:Chebyshev-RQI Subspace}
	\begin{algorithmic}[1]
		\REQUIRE $\matA$, $\matB$, \textit{NEV}, $m$,
                 ${\rm dim}_{\max}$, $\rm It_{\max}$,
                 $\varepsilon_{\text{\tiny{CRS}}}$,
                 $\varepsilon_{\text{\tiny{RQI}}}$,
                 $\rm It_{\text{max-linear}}$,
                 $\vecv_0$.\\
		\ENSURE $\mathbf{\Lambda}$, $\matW$, $\rm It_{\text{total}}$.\\
		\STATE $\vecx = \vecv_0, \rm It_{\text{total}} = 0, \matW=[\ ], \mathbf{\Lambda} = [\ ]$.
	
		\FOR{$n=1$ to \textit{NEV}}
    		\STATE Let $k = 1$, $\matV_k = \{ \vecx \}$,
                   $\theta^{(k)}_n= \frac{\vecx^{\top} \matA \vecx}{\vecx^{\top} \matB \vecx}$,
                   \revise{$\vecr = \frac{(\matA - \theta^{(k)}_n \matB) \vecx}{\theta^{(k)}_n \Vert \vecx \Vert}$.}
    		\STATE Let $a = 0$, $b = 0$ and $\tilde{\sigma}_1 = 0$.
    		\FOR{$k = 1$ to $\rm It_{\max}$}
        		\IF{$ \revise{\Vert \vecr \Vert} < \varepsilon_{\text{\tiny{CRS}}}$}
            		\STATE $\rm It_{\text{total}} = k + \rm It_{\text{total}}$,
                           $\matW = [\matW, \vecx]$, $\mathbf{\Lambda} = [\mathbf{\Lambda}, \theta^{k}_n]$.
            		\STATE break.
        		\ENDIF
        		\IF{$\text{dim}(\matV_k) \geq {\rm dim}_{\max}$}
		        \STATE \texttt{RESTART}: $\matV_k=[ \vecx ]$, $\theta^{(k)}_n = \frac{\vecx^{\top} \matA \vecx}{\vecx^{\top} \matB \vecx}$.
		        \ENDIF
        		\STATE $\matC = \matA - \theta^{(k)}_n \matB$.
        		\IF{$k=1$}
        		    \STATE $\vecz = \matC \vecx$.
        		\ELSE
        		    \STATE $\vecz = {\texttt{ChebyshevFilter}}(\matC, \vecx, m, a, b, \tilde{\sigma}_1)$.
        		    \STATE $\vect = \texttt{IRQI}(\matC, \vecx, \rm It_{\text{max-\tiny{RQI}}} = 1, \varepsilon_{\text{\tiny{RQI}}}, \rm It_{\text{max-linear}})$.
        		\ENDIF
        		\STATE $(\vecz, \vect) = {\tt Orthonormal}(\matW, \vecz, \vect)$.
        		\STATE $(\vecz, \vect) = {\tt Orthonormal}(\matV_k, \vecz, \vect)$.
        		\STATE Let $\matV_{k+1} = \{\matV_k, \vecz, \vect\}$.
        		\STATE Let $\tilde{\matA}_{k+1} = \matV_{k+1}^{\top} \matA \matV_{k+1}$ and $\tilde{\matB}_{k+1} = \matV_{k+1}^{\top} \matB \matV_{k+1}$.
        		\STATE Compute the eigenvalues of the matrix $\tilde{\matA}_{k+1} - \theta^{(k)}_n \tilde{\matB}_{k+1}$
        		       to determine the filter parameters $a$, $b$, and $\tilde{\sigma}_1$.
        		\STATE Compute the eigenpairs $(\mu_i, \vecy_i)$ of the matrix pair
                       $(\tilde{\matA}_{k+1}, \tilde{\matB}_{k+1})$ in ascending  order,
        		       $$
        		       \tilde{\matA}_{k+1} \vecy_i = \mu_i \tilde{\matB}_{k+1} {\vecy_i},
        		       \quad
                       i=1, 2, \ldots, k+1.
                       $$
        		\STATE Update $\vecx = \matV_{k+1} \vecy_1$, $\theta^{(k+1)}_n = \mu_1$.
                \STATE \revise{Update $\vecr = \frac{(\matA - \theta^{(k+1)}_n \matB) \vecx}{\theta^{(k+1)}_n \Vert \vecx \Vert}$.}
    		\ENDFOR
		    \STATE $\matV^{n} = \matV_{k}$, $\vecx = \matV^{n} {\vecy_2}$.
		\ENDFOR
	\end{algorithmic}
\end{algorithm}

In Algorithm~\ref{alg:Chebyshev-RQI Subspace},
the meaning of the notations are the same as \revise{that} in Algorithm 4.
Line~1 is the initialization for the algorithm;
Line~3--4 is the initialization for computing the $n$-th eigenpair.
To compute the first smallest eigenpair, a random initial vector $\vecv_0$ can be used;
To compute the $(n+1)$-th eigenpair, it prefers to
use $\matV^{n} {\vecy_2}$ to start the iteration process.
See Line~29 in the algorithm.

Line~20--26 are the Rayleigh-Ritz process. In this process,
first an orthonormal basis $\matV_{k+1}$ is constructed for the projection subspace;
then the matrix pencil $(\matA, \matB)$ is projected onto the subspace
and the projected matrix pencil $(\matV_{k+1}^{\top} \matA \matV_{k+1}$,
$\matV_{k+1}^{\top} \matB \matV_{k+1})$ is constructed.
\rrevise{We should note that $\text{dim}(\matV_{k+1}) = (2k + 1)\%{\rm dim}_{\max}$, and in our experiments we set ${\rm dim}_{\max}=80$, which means the dimension of projected dense eigenproblem should be small in the iteration. The computational cost would be no significant increase in Rayleigh-Ritz process, although we add two augmented vectors into subspace at a time.
Besides, due to the much faster convergence speed than that of CD, the total number of Rayleigh-Ritz process will decrease significantly in CRS method, please refer to Tables~\ref{tab:scale1}--\ref{tab:scale3}.
}
And finally, \revise{in Line~26} the product $\matV_{k+1} \vecy_1$ is used to approximate the desired eigenvector.
In Line~24, the interval selection strategy is the same as Chebyshev Davidson method discussed in Sect.~\ref{sect:Cheby-davidson}, i.e., we let $\tilde{\sigma}_{1}$, $a$, and $b$ represent the smallest,
the second smallest and the largest eigenvalues of the matrix $\matV_{k+1}^{\top} \matC \matV_{k+1}$.

It is necessary to point out that we don't have to explicitly compute the new projected matrix pencil
$(\tilde{\matA}_{k+1}, \tilde{\matB}_{k+1})$ in Line~23.
Actually, we note that
\begin{eqnarray}
\label{eqn:def-tilde-A}
\tilde{\matA}_{k+1} =
\matV_{k+1}^{\top} \matA \matV_{k+1} =
\left[
\begin{array}{ccc}
  \matV_k^{\top} \matA \matV_k & \matV_k^{\top} \matA \vecz & \matV_k^{\top} \matA \vect \vspace{3pt} \\
  \vecz^{\top} \matA \matV_k & \vecz^{\top} \matA \vecz & \vecz^{\top} \matA \vect \vspace{3pt} \\
  \vect^{\top} \matA \matV_k & \vect^{\top} \matA \vecz & \vect^{\top} \matA \vect
\end{array}
\right],
\end{eqnarray}
and
\begin{eqnarray}
\label{eqn:def-tilde-B}
\tilde{\matB}_{k+1} =
\matV_{k+1}^{\top} \matB \matV_{k+1}=
\left[
\begin{array}{ccc}
  \matV_k^{\top} \matB \matV_k & \matV_k^{\top} \matB \vecz & \matV_k^{\top} \matB \vect \vspace{3pt} \\
  \vecz^{\top} \matB \matV_k & \vecz^{\top} \matB \vecz & \vecz^{\top} \matB \vect \vspace{3pt} \\
  \vect^{\top} \matB \matV_k & \vect^{\top} \matB \vecz & \vect^{\top} \matB \vect
\end{array}
\right].
\end{eqnarray}
Since $\matV_k^{\top} \matA \matV_k$ and $\matV_k^{\top} \matB \matV_k$ in~\eqref{eqn:def-tilde-A} and~\eqref{eqn:def-tilde-B}
already exist after the previous iteration,
also note that the matrices $\tilde{\matA}_{k+1}$ and $\tilde{\matB}_{k+1}$ are symmetric,
we only need to compute the diagonal blocks at $(2,2)$ and $(3,3)$,
and upper triangular parts,
i.e., we only need to compute
four vectors
$\matV_k^{\top} \matA \vecz$, $\matV_k^{\top} \matA \vect$,
$\matV_k^{\top} \matB \vecz$, and $\matV_k^{\top} \matB \vect$,
six scalars
$\vecz^{\top} \matA \vecz$,
$\vecz^{\top} \matA \vect$,
$\vect^{\top} \matA \vect$,
$\vecz^{\top} \matB \vecz$, $\vecz^{\top} \matB \vect$,
and $\vect^{\top} \matB \vect$.

One crucial step of the algorithm is Line~18, where an augmentation vector $\vect^{(k)}$
is obtained by calling one step Rayleigh quotient iteration,
and the linear equations in Rayleigh quotient iteration is solved inexactly by
an iterative method. See Line~4 in Algorithm~\ref{alg:IRQI}.

It is worth mentioning that we tried different combinations of the various Krylov subspace methods and preconditioners in PETSc library~\cite{petsc} for solving the Rayleigh quotient equation.
Finally we find that CRM (Conjugate Residual Method) and MINRES (Minimal Residual Method)
without preconditioning show good results, and CRM performs better than MINRES.
Simoncini~\cite{simoncini2002inexact} and Jia~\cite{jia2012convergence}
have made some convergence analysis on the application of these two methods
in inexact Rayleigh quotient iteration.
However, to make the algorithms work efficiently,
either a good approximate eigenpair should be provided,
or the shifted matrix should be decomposed to construct a preconditioner.
In real applications, it is usually difficult to give a good
approximate eigenpair, and for solving large scale problem,
it is too expensive to decompose the shifted matrix.
Furthermore, they only considered to compute one
eigenvalue and \revise{associated} eigenvector.

In our method, we only control \revise{the} iteration number of the Krylov method by setting the maximal number of iterations as $\rm It_{\text{max-linear}}$.
For solving the RQI equation inexactly, the initial guess of inner iteration will have strong influence on the performance of the method, and the zero vector would be a good choice.

There are two clear characters for CRS algorithm.
Firstly, compared with Chebyshev-Davidson method,
the filter interval for CRS can be selected more flexible because the error in the eigendirection
corresponding to the larger eigenvalue can be easily eliminated by using RQI process.
The second characteristic is that Chebyshev iteration and RQI are complementary to each other
in CRS method. The Chebyshev iteration could always provide a good start vector for RQI,
and RQI could obtain more approximate information about \revise{the} 
 desired eigendirection,
which in turn helps the algorithm to obtain a better filter interval for next iteration step.
Therefore, CRS algorithm is more flexible and stable
than the Chebyshev-Davidson algorithm,
and it should converge faster than Chebyshev-Davidson algorithm.

\section{Numerical Experiments}
\label{sect:numer-result}

In this section, we present some numerical results to compare Chebyshev-RQI subspace (CRS) method with some other eigenvalue solution methods for computing several smallest eigenpairs
of the symmetric GEP~(\ref{eq:GEP}).
The compared methods include Chebyshev-Davidson (CD) method~\cite{miao2020chebyshev},
Krylov-Schur (KrylovSchur) method~\cite{stewart2002krylov,slepc},
the efficient Jacobi-Davidson (JD) method~\cite{sleijpen2000jacobi,slepc},
the locally optimal block preconditioned conjugate gradient (LOBPCG) method~\cite{knyazev2001toward,slepc} and
the Generalized Davidson (GD) method~\cite{morgan1986generalizations,slepc}.
KrylovSchur, JD and LOBPCG are block type methods,
i.e, they can compute several eigenpairs simultaneously
while other methods only compute one eigenpair each time.
We will show the performance of the aforementioned six methods
with respect to the outer iteration steps (\textbf{IT}),
the number of matrix-vector products (\textbf{MV}) and
the computing time in seconds (\textbf{TIME}).

In the following, $\vecr^{(k)} = \frac{(\matA - \theta^{(k)} \matB) \vecx^{(k)}}{\theta^{(k)} \Vert \vecx^{(k)} \Vert}$
represents the relative residual vector of CRS algorithm with $\theta^{(k)}$
being the Rayleigh quotient associated with the $k$-th iterate approximate eigenvector $\vecx^{(k)}$.
And the whole iteration process will be terminated as long as their current relative residual
norm is less than the prescribed stopping criterion 
\revise{$\| \vecr^{(k)}\Vert < 10^{-10}$}.

\subsection{The eigenvalue problem}
Consider the following two-dimension beam free vibration system
\begin{eqnarray}
\label{eq:governing equation}
\left\{
\begin{array}{rcl}
\bm{\sigma} \cdot \nabla -\rho \ddot{\mathbf{u}} &=& 0, \quad \text{in} \; \Omega, \\
\bm{\sigma} \cdot \mathbf{n} = \mathbf{\overline{t}}  &=& 0, \quad \text{on} \; \Gamma_N, \\
\mathbf{u}                                       &=& 0, \quad \text{on} \; \Gamma_D,
\end{array}
\right.
\end{eqnarray}
where $\mathbf{u}$ is the displacement, 
$\bm{\sigma}$ is the stress, 
$\rho$ is the density of material.
The computing domain is $\Omega = [0,10] \times [0,2]$, 
$\Gamma_D = \left\{ (x,y) \vert x=0,\ y \in [0,2] \right\}$ and 
$\Gamma_N = \partial\Omega \setminus \Gamma_D$, 
$\mathbf{n}$ is the outer normal vector of $\partial \Omega$.

Let
\begin{eqnarray*}
\mathbf{u} = \begin{bmatrix}
u(x,y,t)	\\
v(x,y,t)
\end{bmatrix},
\qquad
\delta \mathbf{u} = \begin{bmatrix}
	w(x,y,t)	\\
	s(x,y,t)
	\end{bmatrix},	
\end{eqnarray*}
denote the displacement and the virtual displacement, respectively.
By the principle of virtual work, we have
\begin{eqnarray*}
&& -\int_\Omega (\bm{\sigma} \cdot \nabla - \rho \ddot{\mathbf{u}}) \delta \mathbf{u} \df \Omega + \int_{\Gamma_N}(\bm{\sigma} \cdot \mathbf{n} -\mathbf{\overline{t}})  \delta \mathbf{u} \df s\\
=&& -\int_\Omega (\bm{\sigma} \cdot \nabla -\rho \ddot{\mathbf{u}}) \delta \mathbf{u} \df \Omega + \int_{\Gamma_N}(\bm{\sigma} \cdot \mathbf{n} - \mathbf{\overline{t}}) \delta \mathbf{u} \df s
+ \int_{\Gamma_D}(\bm{\sigma} \cdot \mathbf{n} ) \delta \mathbf{u}  \df s\\
=&&-\int_\Omega  (\bm{\sigma} \cdot \nabla )\delta \mathbf{u} \df \Omega + \int_{\partial\Omega}(\bm{\sigma} \cdot \mathbf{n} ) \delta \mathbf{u} \df s
+ \int_\Omega \rho \ddot{\mathbf{u}} \delta \mathbf{u} \, \df \Omega -  \int_{\Gamma_N} \mathbf{\overline{t}} \delta \mathbf{u} \df s\\
=&& \int_\Omega \bm{\sigma} \colon\frac{1}{2} \left (\delta \mathbf{u} \nabla + \delta \mathbf{u}^{\top} \nabla \right) \df \Omega + \int_\Omega \rho \ddot{\mathbf{u}}\delta \mathbf{u} \df \Omega
-  \int_{\Gamma_N} \mathbf{\overline{t}} \delta \mathbf{u} \df s \\
=&& \int_\Omega 2 \mu \bm{\varepsilon}(\mathbf{u}) \colon \bm{\varepsilon}(\delta \mathbf{u}) + \lambda \nabla \cdot \mathbf{u} \nabla \cdot \delta \mathbf{u}  \df \Omega +\int_\Omega \rho \ddot{\mathbf{u}} \delta \mathbf{u} \df \Omega=0,
\end{eqnarray*}
where $\bm{\varepsilon}$ is the strain of the material.
We should note that,
\begin{eqnarray*}
\int_\Omega \bm{\sigma} \colon \frac{1}{2} \left (\delta \mathbf{u} \nabla + \delta \mathbf{u}^{\top} \nabla \right ) \df \Omega 
= -\int_\Omega  (\bm{\sigma} \cdot \nabla ) \delta \mathbf{u} \df \Omega + \int_{\Gamma}^{}(\bm{\sigma} \cdot \mathbf{n} ) \delta \mathbf{u} \df s.
\end{eqnarray*}
After discretization by finite element method, we have
discrete equation for (\ref{eq:governing equation})
\begin{eqnarray*}
\mathbf{M} \ddot{\mathbf{u}}_h(t) + \mathbf{K} \mathbf{u}_h(t) = 0,
\end{eqnarray*}
where $\mathbf{u}_h$ is the displacement of nodal points,
$\mathbf{K}$ and $\mathbf{M}$ are stiff matrix and mass matrix,
which \revise{are} defined by variational form as follows
\begin{eqnarray*}
&& \mathbf{K} = 
\int_\Omega \left\{ 2 \mu \left[ \frac{\partial u}{\partial x} \cdot           
             \frac{\partial w}{\partial x} + \frac{\partial v}{\partial y} \cdot 
             \frac{\partial s}{\partial y} + 
             \frac{1}{2} \left( \frac{\partial v}{\partial x} + \frac{\partial u}{\partial y} \right) \cdot \left( \frac{\partial s}{\partial x} + \frac{\partial w}{\partial y} \right) \right] \right.\\
 && \qquad\quad  + \left. \lambda \left( \frac{\partial u}{\partial x} + 
             \frac{\partial v}{\partial y} \right) \cdot \left( 
             \frac{\partial w}{\partial x} + 
             \frac{\partial s}{\partial y} \right) \right\} \df \Omega,  \\
&&\mathbf{M} = \int_\Omega \rho ( u \cdot w + v \cdot s) \df \Omega,
\end{eqnarray*}
where the  Lam\'{e} constants $\mu$ and $\lambda$ can be computed by the Young's modulus $E$ and Poisson's ratio $\nu$ of the material, i.e.,
$\mu = \frac{E}{2(1+\nu)}$ and $\lambda = \frac{E \nu}{(1+\nu)(1 - 2\nu)}$.

By using {\tt FreeFem++} software \cite{FreeFem}, we generate the stiffness matrix $\mathbf{K}$
and the mass matrix $\mathbf{M}$ that correspond to the matrices $\matA$ and $\matB$ in the symmetric GEP (\ref{eq:GEP}), respectively.
In the following numerical results, three scales of the problem will be tested:
S1: $N = 46958$;
S2: $N = 187778$;
S3: $N = 1143146$;
where $N$ is the number of the DOFs.

\begin{figure}[htbp]
\centering
\includegraphics[scale=0.3]{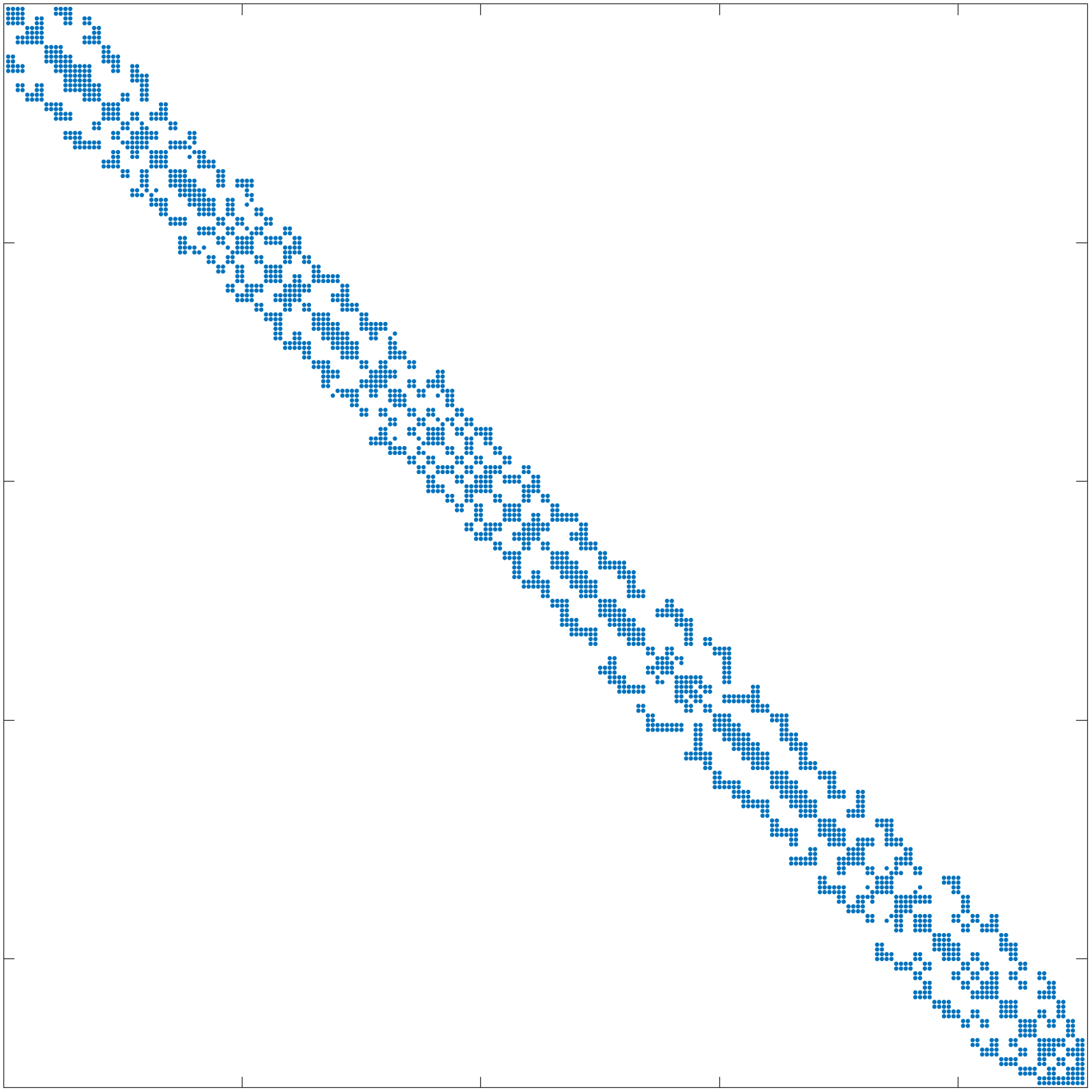}
\hspace{3em}
\includegraphics[scale=0.3]{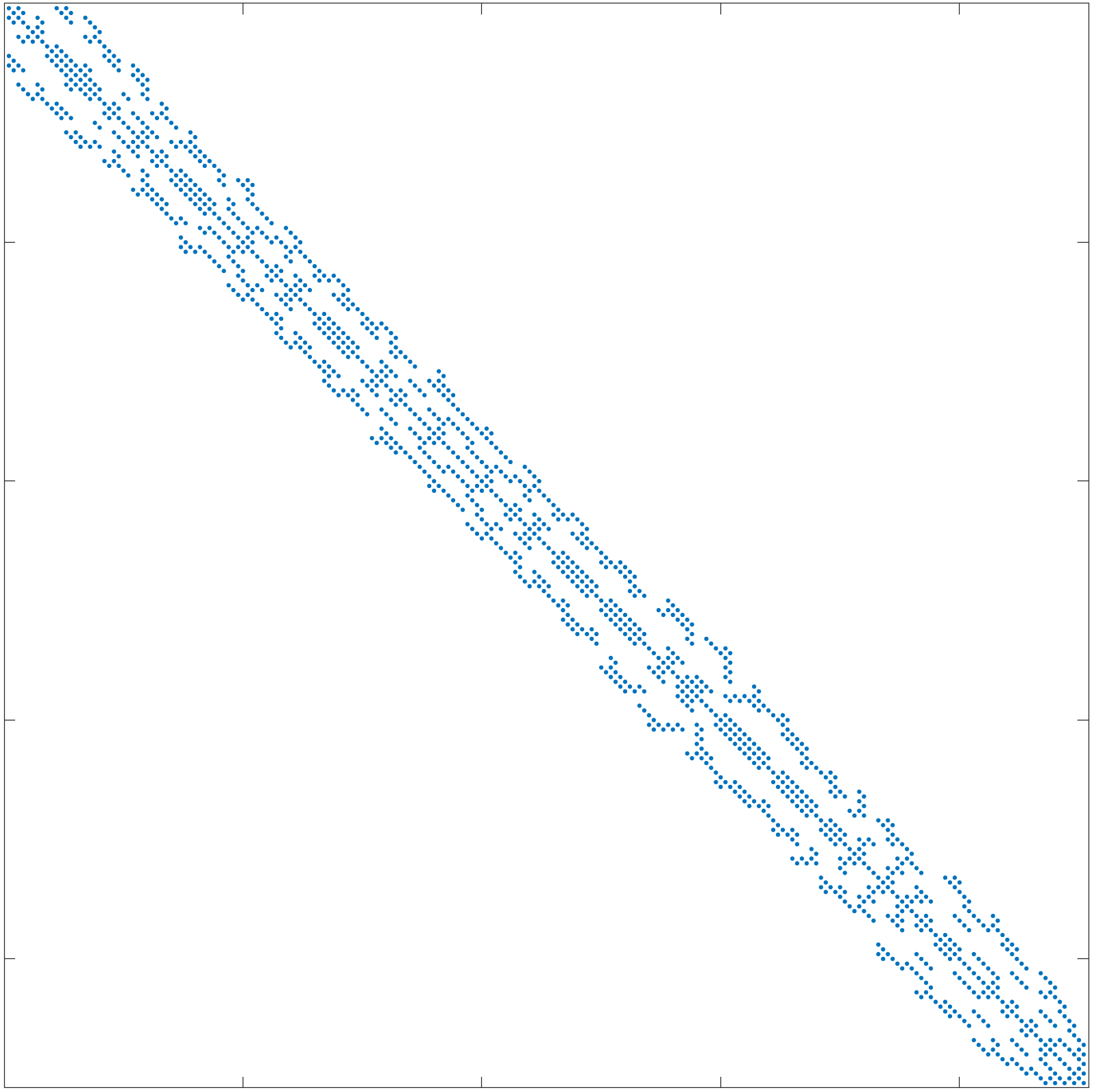}
\caption{The structure of matrix pair $A$ (left) and $B$ (right) for $N=187778$.}
\end{figure}

\subsection{Numerical Results}

The numerical experiments are carried out by using
PETSc~\cite{petsc} and SLEPc~\cite{slepc} in which
MPI based distributed vectors and sparse matrices are provided.
Besides CD and CRS method, numerical experiments for other compared methods
are carried out by calling SLEPc.
We set \texttt{-eps\_ncv 120} when $NEV = 20$,
and \texttt{-eps\_ncv 200} when $NEV = 100$ for KrylovSchur method
because this method performs worst or even fail with default setting,
here \texttt{-eps\_ncv} denotes the largest dimension of worrking subspace in KrylovShur method.
The default parameters \rrevise{which are provided by SLEPc} are used for all \revise{other cases}. 

The CRM (Conjugate Residual Method) without preconditioner is used as the linear solver, and zero vector is used as the initial vector for solving the RQI equation~(\ref{single-step IRQI}) in CRS method. We use the mark ``-'' to indicate that the computing time is more than two hours.

Numerical experiments were carried out on a cluster, all blade nodes are equipped with two 2.60GHz Intel(R) Xeon(R) Gold 6132 CPU, each CPU has 14 cores, and each node has 28 cores with 12 $\times$ 8GB DDR4 2400MHz ECC total 96GB memory.

\begin{table}[htbp]
	\centering
	\caption{Numerical results for S1.}
	  \begin{tabular}{cccccccc}
	  \toprule
	  N     & METHOD & \multicolumn{3}{c}{NEV=20} & \multicolumn{3}{c}{NEV=100} \\
  \cmidrule{3-8}          &       & \textbf{IT}    & \textbf{MV}    & \textbf{TIME(s)}  & \textbf{IT}    & \textbf{MV}    & \textbf{TIME(s)} \\
	  \midrule
	  \multirow{6}[2]{*}{46958} & KrylovSchur & 153  & 26980 & 371.94 & 162   & 41565 & 857.92 \\
			& LOBPCG & 125   & 14837 & 85.05    & 526   & 55806 & 336.00 \\
			& GD    & 4309  & 12375 & 83.76 & 15166 & 45103 & 996.25 \\
			& JD    & 136   & 40736 & 109.06 & 617   & 187075 & 609.32 \\
			& CD   & 1593  & 54272 & 114.76 & 6815  & 232242 & 543.66 \\
			& CRS   & 408   & 22167 & \textbf{63.74 }& 1851  & 146996 & \textbf{309.86} \\
	  \bottomrule
	  \end{tabular}%
	\label{tab:scale1}%
  \end{table}%

In Table~\ref{tab:scale1}, we report the numerical results of scale S1 
for computing $NEV$ smallest eigenvalues and the
corresponding eigenvectors for the symmetric generalized eigenvalue problem \eqref{eq:GEP}
by using 1 processor. The restart number and the polynomial order involved in the
CD and CRS method are set to be 80 and 30, respectively.
And the number of iteration for solving RQI linear equation in CRS is set to 50.

From Table~\ref{tab:scale1}, 
we observe that KrylovSchur method takes the most time when $NEV=20$ and LOBPCG is more effective than JD.
The CRS method is nearly twice faster than CD in term of computing time, and four times faster in term of number of iterations, respectively.
JD and LOBPCG are slower than GD when $NEV=20$, but much faster when $NEV=100$. This shows that the block algorithm
has great advantages for computing multiple eigenpairs.
And we observe that among the seven methods the CRS method performs the best
for both cases of $NEV=20$ and $NEV=100$,
even though CRS is not a block algorithm.

\begin{table}[htbp]
	\centering
	\caption{Numerical results for S2.}
	  \begin{tabular}{cccccccc}
	  \toprule
	  N     & METHOD & \multicolumn{3}{c}{NEV=20} & \multicolumn{3}{c}{NEV=100} \\
  \cmidrule{3-8}          &       & \textbf{IT}    & \textbf{MV}    & \textbf{TIME(s)}  & \textbf{IT}    & \textbf{MV}    & \textbf{TIME(s)} \\
	  \midrule
	  \multirow{6}[2]{*}{187778} & KrylovSchur & 379     & 66556   & 1630.36   & 357   & 89844 & 2258.28 \\
			& LOBPCG & 332   & 38125 & 58.92 & 1421  & 138448 & 223.09 \\
			& GD    & 9812  & 28170 & 54.49 & 35450 & 105523 & 653.66 \\
			& JD    & 152   & 46176 & 36.46 & 697   & 215653 & 196.98 \\
			& CD   & 2501  & 110598 & 59.46 & 10558 & 466763 & 269.67 \\
			& CRS   & 456   & 40804 & \textbf{27.03} & 2026  & 251690 & \textbf{125.90} \\
	  \bottomrule
	  \end{tabular}%
	\label{tab:scale2}%
  \end{table}%

In Table~\ref{tab:scale2}, we report the numerical results of scale S2 for computing $NEV$ smallest eigenvalues and the
corresponding eigenvectors for the symmetric generalized eigenvalue problem (\ref{eq:GEP}) by using 18 processors. The restart number and the polynomial order involved in the
CD and CRS method are also set to be 80 and 40, respectively.
The number of iteration for solving RQI linear equation in CRS is set to 90.

From Table~\ref{tab:scale2}, one can see that in both $NEV=20$ and $NEV=100$ cases,
CRS performs the best by comparing the solution time.
KrylovSchur is much slower than the other five methods.
For this scale problem, JD method is more efficient than LOBPCG method.
The number of iterations and computing time of CD are
about 5 times and 2.2 times of CRS respectively.

\begin{table}[htbp]
	\centering
	\caption{Numerical results for S3.}
	  \begin{tabular}{cccccccc}
	  \toprule
	  N     & METHOD & \multicolumn{3}{c}{NEV=20} & \multicolumn{3}{c}{NEV=100} \\
  \cmidrule{3-8}          &       & \textbf{IT}    & \textbf{MV}    & \textbf{TIME(s)}  & \textbf{IT}    & \textbf{MV}    & \textbf{TIME(s)} \\
	  \midrule
	  \multirow{6}[2]{*}{1143146} & KrylovSchur & -     & -     & -     & -     & -     & - \\
			& LOBPCG & 1062  & 128744 & 258.37 & 4820  & 465745 & 1026.06 \\
			& GD    & 25275 & 72567 & 195.19 & 93338 & 277951 & 2148.82 \\
			& JD    & 212   & 67219 & \textbf{79.27} & 794   & 250761 & 336.71 \\
			& CD   & 6326  & 281991 & 202.01 & 27676 & 1326207 & 893.43 \\
			& CRS   & 500   & 116908 & 83.57 & 2062  & 545207 & \textbf{332.87} \\
	  \bottomrule
	  \end{tabular}%
	\label{tab:scale3}%
  \end{table}%

In Table~\ref{tab:scale3}, we report the numerical results of scale S3 
for computing $NEV$ smallest eigenvalues and the corresponding eigenvectors
for the symmetric generalized eigenvalue problem~\eqref{eq:GEP}
by using 112 processors. The restart number and the polynomial order involved in the
CD and CRS method are also set to be 80 and 40, respectively.
And the number of iteration for solving RQI linear equation in CRS is set to 250.
We should remark that, in order to minimize the impact of communication between nodes,
all experiments are done on the same 4 nodes.

From Table~\ref{tab:scale3}, we observe that KrylovSchur method fail.
In this case, JD method is much more efficient,
and it is about 3 times faster than LOBPCG method in term of computing time.
The computing time of CD method are about 2 times more than that of CRS.
JD and CRS methods perform best in this case.	
By comparing the results of CD and CRS from Tables~\ref{tab:scale1}--\ref{tab:scale3},
we can find that CRS performs better and better than CD with the increase of problem size.

\begin{figure}[htbp]
	\centering
	\includegraphics[scale=0.3]{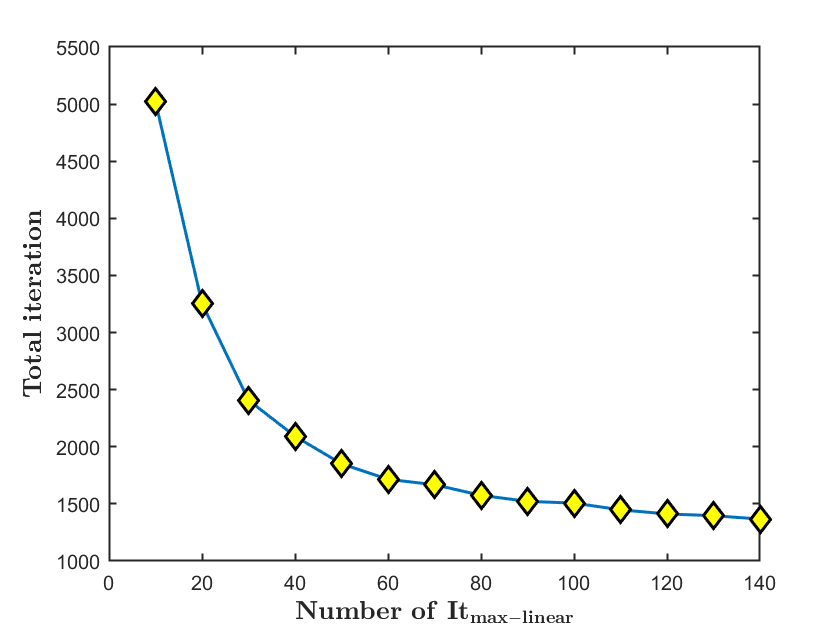}
	\hspace{0.1in}
	\includegraphics[scale=0.3]{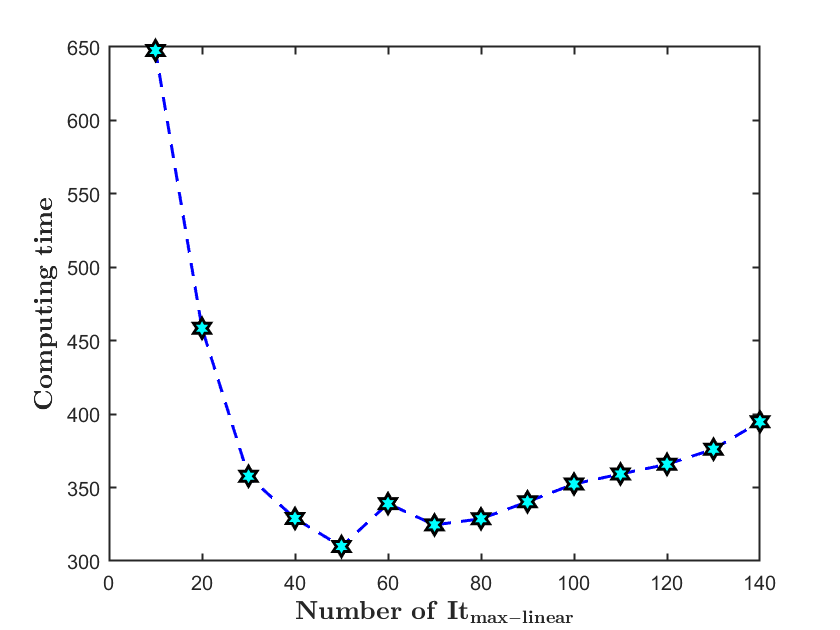}
	\caption{The total iteration (left) and computing time (right) of CRS method with different $\rm It_{\text{max-linear}}$ for S1 ($N=46958$, 1 processor).}
	\label{inner_refine4}
\end{figure}
\begin{figure}[htbp]
	\centering
	\includegraphics[scale=0.3]{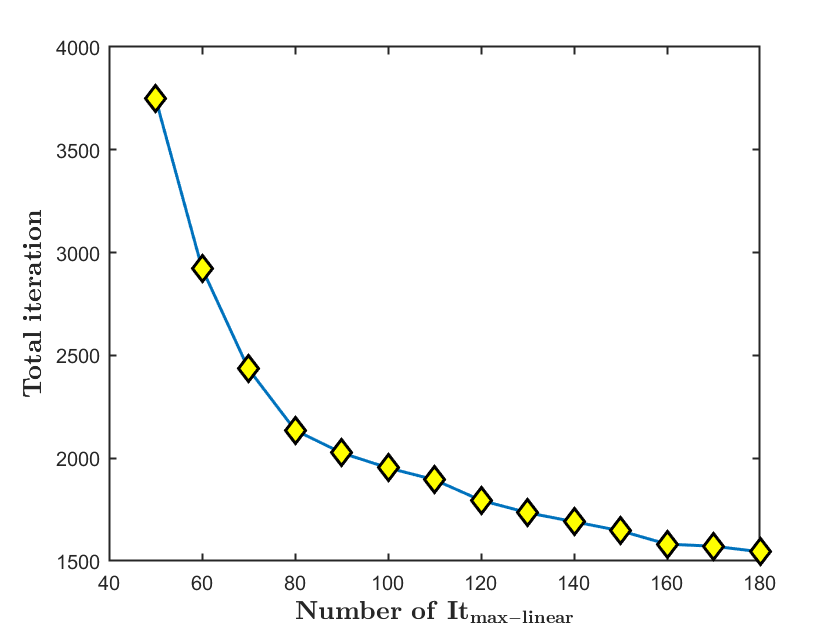}
	\hspace{0.1in}
	\includegraphics[scale=0.3]{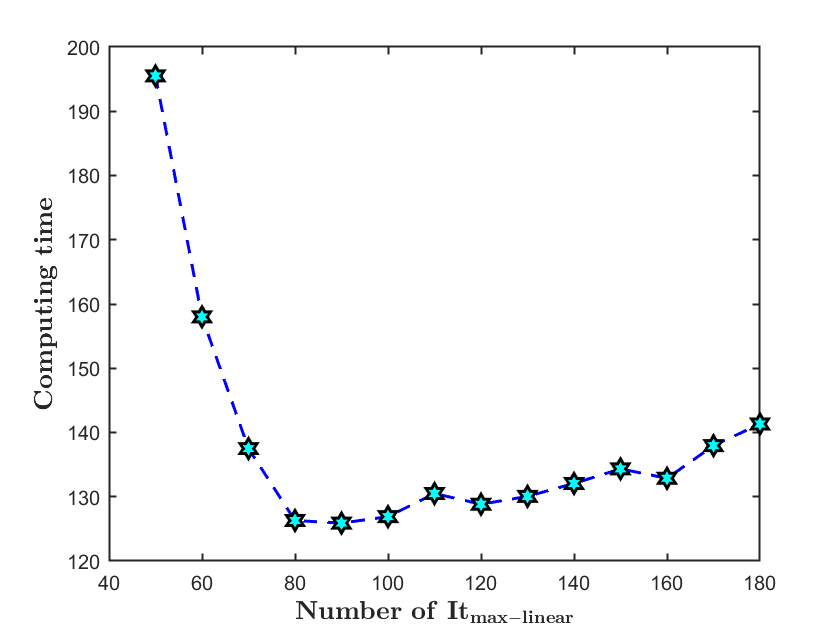}
	\caption{The total iteration (left) and computing time (right) of CRS method with different $\rm It_{\text{max-linear}}$ for S2 ($N=187778$, 18 processors).}
	\label{inner_refine5}
\end{figure}
\begin{figure}[!htbp]
	\centering
	\includegraphics[scale=0.3]{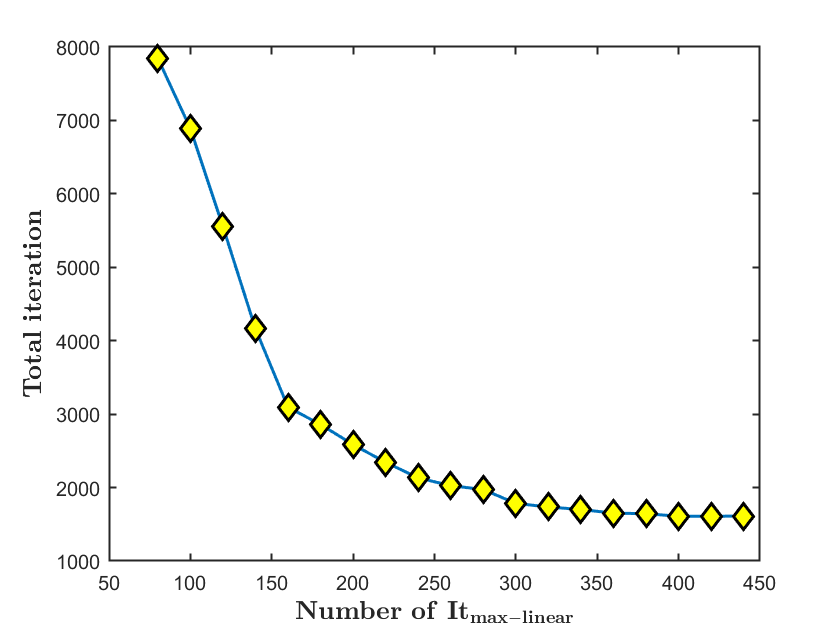}
	\hspace{0.1in}
	\includegraphics[scale=0.3]{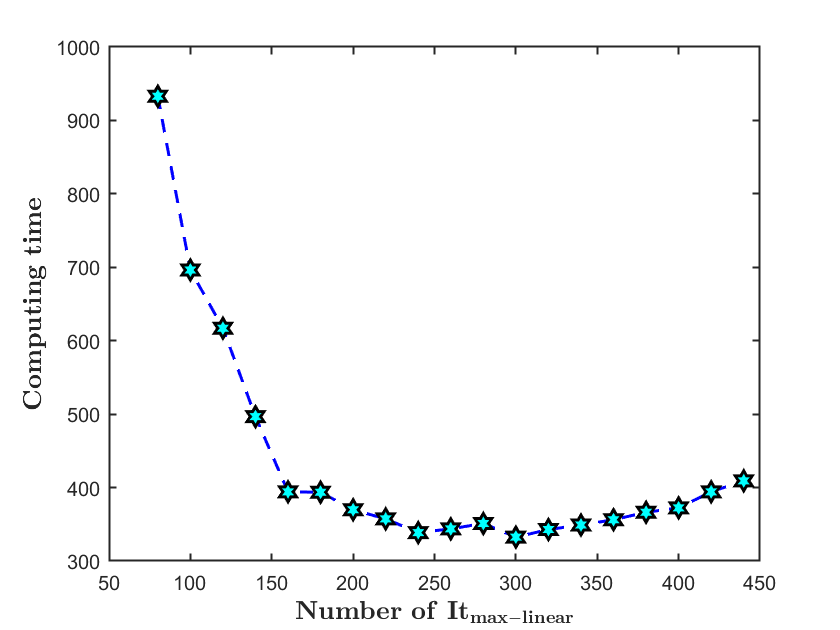}
	\caption{The total iteration (left) and computing time (right) of CRS method with different $\rm It_{\text{max-linear}}$ for S3 ($N=1143146$, 112 processors).}
	\label{inner_refine6}
\end{figure}

The number of linear iteration for solving the RQI equation
is the most important ingredient
affecting the performance of CRS method.
In Figures~\ref{inner_refine4}, \ref{inner_refine5}, and \ref{inner_refine6},
the iteration numbers and computing time for solving \revise{these}  three scale problems
are plotted with the increase of $\rm It_{\text{max-linear}}$,
the maximal iteration number for solving the RQI linear equation.
From these figures, we can find that, when $NEV=100$, 
the optimal value of $\rm It_{\text{max-linear}}$ is about 
50 for S1 scale, 90 for S2 scale, and 250 for S3 scale
in terms of computing time 
(where the polynomial order $m=30$ for S1, and $m=40$ for S2 and S3).
The data in Tables~\ref{tab:scale1}--\ref{tab:scale3} are obtained
by using these optimal values.
From these figures, one can also find that when the iteration number
$\rm It_{\text{max-linear}}$ is greater than the optimal value, 
the CRS iteration decreases very little and the computing time increases not too much.

\begin{figure}[htbp]
	\centering
	\includegraphics[scale=0.33]{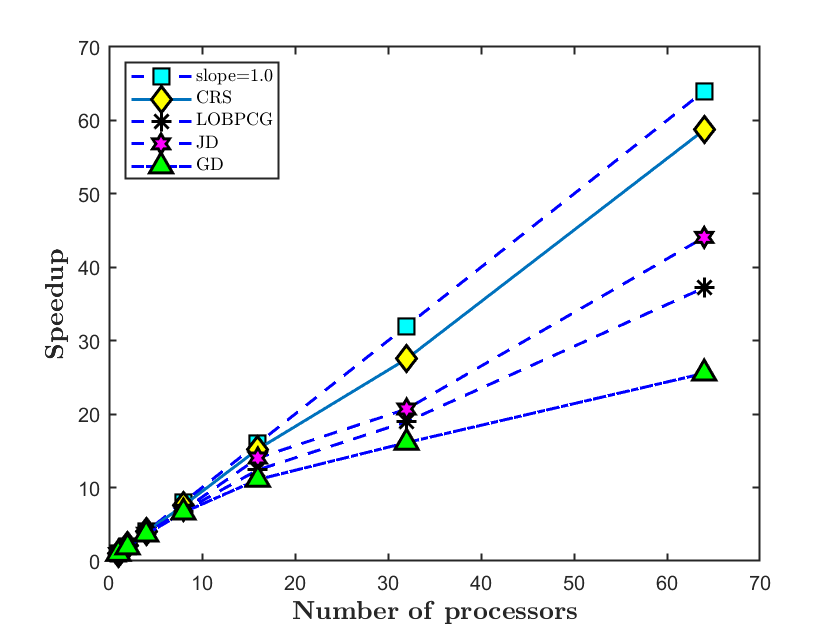}
	\hspace{0.1in}
	\includegraphics[scale=0.33]{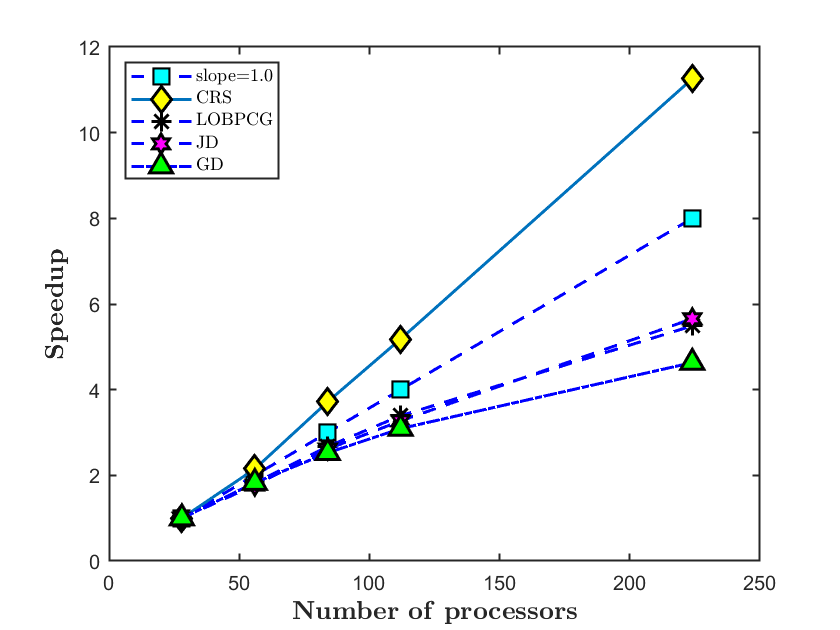}
	\caption{Strong scalability  for four methods (left: $N=187778$; right: $N=1143146$).}
	\label{strongscalability}
\end{figure}

In the following, we show the strong scalability of several methods mentioned above \revise{for computing} the first 20 smallest eigenpairs for two scales: $N=187778$ and $N=1143146$.
For these two scale cases, the linear iteration number for solving the RQI equations
is set to 90, and 250, respectively.
The speedup curves for the strong scalability are plotted in Fig.~\ref{strongscalability}.
From this figure we observe that, among the four methods CRS has the best strong scalability.
The line of CRS is almost linear when $N=187778$ (left) and
superlinear when $N=1143146$ (right).
We analyzed the time of each module in the algorithm and found that
the time cost of {\texttt{ChebyshevFilter}} and CRM iteration decrease superlinearly.

At last in this section, it is necessary to point out that the number of restart is not
a sensitive parameter in CRS. In numerical experiments, the algorithm converges within
40 iteration steps for computing one of the eigenpairs, i.e.,
the dimension of projection subspace would not exceed 80 and
there is no need to restart.
The restart parameter is set to prevent excessive storage for computing \revise{the eigenpairs which converge very slowly}. 
\section{Concluding Remarks}
\label{sect:conclusion}

In this paper, a new subspace algorithm for solving symmetric generalized eigenvalue problems
is obtained by combining the technique of Chebyshev polynomial filter
and inexact Rayleigh quotient iteration in the iteration process.
The obtained method is named as CRS algorithm that can be used to
compute several smallest eigenvalues and the corresponding eigenvectors.

Numerical results for a kind of vibration model
show that the performance of the proposed algorithm is \revise{very good}.
Compared with some eigenvalue algorithms implemented in SLEPc,
CRS method can compute multiple eigenpairs faster,
and also it shows the best parallel scalability in our experiments.
Therefore, CRS method has the potential to solve large scale eigenvalue problem
in modal analysis of mechanical vibration.

In CRS method, it is needed to solve a linear equation related to RQI
in each CRS iteration,
and this is the dominant cost of CRS algorithm.
For solving the linear equations, a Krylov subspace method (CRM) is employed.
We studied the influence of the iteration number of Krylov method on the performance
of the whole CRS algorithm.
Numerical results show that
a relatively optimal maximal iteration number is related to the scale of the problem.

For CRS algorithm, the following issues should be studied further:
\begin{itemize}
\item The method CRS proposed in this paper is a kind of one vector algorithm,
      that is, only one eigenpair can be obtained in each iteration.
      To compute multiple eigenpairs, block type method
      is usually preferable. It is deserved to
      consider the generalization of CRS to block version \revise{and the corresponding efficient implementation techniques} for computing several eigenvalues simultaneously.
\item In CRS method, one dominant cost is to solve the linear equation
      of RQI process. To trade off the convergence and efficiency,
      the problem is how to design a strategy to adjust
      the number of Krylov iteration adaptively for solving the linear equation.
      By increasing the number of Krylov iterations may
      make the convergence quickly, but at the same time,
      this will result in high computational costs and more computing time.
      This deserves to be studied furthermore to make the algorithm more effective.
\item In this paper, only one iteration step of RQI is implemented in CRS algorithm to       construct the augmented vector.
      One could also get an augmented vector by using several iteration steps of RQI iterations.
      This may result in a better augmented vector and CRS method may converge faster.
      As one linear equation should be solved in each iteration step of RQI,
      so if we use many iteration steps of RQI to produce the augmented vector,
      the cost will be much more expensive.
      Considering the whole efficiency of CRS method,
      this need to be studied further to trade off the iteration steps of RQI and
      the convergence rate of CRS algorithm.
\end{itemize}



\section*{Conflict of interest}

The authors declare that they have no conflict of interest.

\bibliographystyle{spmpsci}      


\bibliography{ChebyshevRQIMethod}

\end{document}